\documentclass[11pt]{article}
\usepackage{a4wide}
\usepackage{amsmath}
\usepackage{amsfonts}
\usepackage{amssymb}
\usepackage{euscript}
\def\boxit{$\sqcap\kern-8pt\sqcup$}
\def\littbox{\hfill{\boxit{}}\vskip .2cm}
\newenvironment{proof}{\noindent {\it Proof.~~}\ }{\  \rule{1mm}{2mm}\medskip}

\newenvironment{proof*}{\noindent {\it Proof.~~}\ }{}
\newtheorem{teo}{Theorem}[section]
\newtheorem{lem}[teo]{Lemma}
\newtheorem{cor}[teo]{Corollary}
\newtheorem{pro}[teo]{Proposition}

\newtheorem{rem}[teo]{Remark}
\newtheorem{ques}[teo]{Question}
\newcommand{\rere}{{\rm I}\!{\rm R}}

\input epsf
\begin{document}
\title{On Tiling Rectangles via the Frobenius Number}

\author { J.L. Ram\'{\i}rez  Alfons\'{\i}n, \\
{\it Universit\'e Pierre et Marie Curie, Paris 6,}\\
{\it Combinatoire et Optimisation} \\
{\it Case 189 - 4 Place Jussieu Paris 75252 Cedex 05, France}}

\maketitle

\begin{abstract} In this paper, we give some sufficient conditions for a $n$-dimensional rectangle to be tiled with a set of bricks.
These conditions are obtained by using the so-called Frobenius number. 
\end{abstract}

\section{Introduction }

Let $a_1,\dots ,a_n$ be positive integers. We denote by $A=(a_1\times\cdots\times a_n)$ the $n$-dimensional rectangle of sides $a_i$, that is, 
$A=\{(x_1,\dots ,x_n)\in\rere^n | 0\le x_i\le a_i, i=1,\dots ,n\}$.
A $n$-dimensional rectangle $A$ is said to be {\em tiled} with {\em bricks} (i.e., small $n$-dimensional rectangles) $B_1,\dots ,B_k$  
if $A$ can be filled entirely with copies of $B_i$, $1\le i\le k$ (rotations allowed). 
It is known \cite{deB,Kla} that rectangle $A=(a_1\times a_2)$ can be tiled with  $(x_1\times x_2)$ if and 
only if $x_1$ divides $a_1$ or $a_2$, $x_2$ 
divides $a_1$ or $a_2$ and if $x_1x_2$ divides one side of $A$ then the other side can be expressed as a nonnegative integer combination of $x_1$ and $x_2$. 
In 1995, Fricke \cite{Fri} gave the following characterization when $n=2$ (see also \cite{BW} for a $n$-dimensional generalization with $k=2$). 

\begin{teo}\label{kler}\cite{Fri} Let $a_1,a_2,x,y$ be positive integers with $g.c.d.(x,y)=1$.  
Then, $(a_1\times a_2)$ can be tiled with $(x\times x)$ and $(y\times y)$ if and only if
either $a_1$ and $a_2$ are both multiple of $x$ or  $a_1$ and $a_2$ are both multiple of $y$ or  
one of the numbers $a_1, a_2$ is a multiple of both $x$ and $y$ and the other can be expressed as a nonnegative integer combination of $x$ and $y$.
\end{teo}

Let us consider the following natural question. 

\begin{ques}\label{qq} Does there exist a function $C=C(x_1,x_2,y_1,y_2)$ such that if $a_1,a_2\ge C$ then 
$(a_1\times a_2)$ can be tiled with $(x_1\times x_2)$ and $(y_1\times y_2)$ for some positive integers $x_1,x_2,y_1$ and $y_2$? 
\end{ques}

An  algebraic result due to Barnes \cite{Bar} seems to show the existence of such $C$. However, Barnes'method does not give an 
explicit lower bound for $C$.

The special case when $x_1=4, x_2=6, y_1=5$ and $y_2=7$ was posed in the 1991 William Mowell Putnam Examination (Problem B-3). 
In this case, Klosinski {\it et. al.} \cite{Klos} gave a lower bound of $C$. Their method was based on knowledge of the {\em Frobenius number}. The {\em Frobenius number}, denoted by $g(s_1,\dots ,s_n)$, of a set of relatively prime positive integers $s_1,\dots ,s_n$, is  defined as the largest integer that is not representable as a nonnegative integer combination of $s_1,\dots ,s_n$. It is well known that

\begin{equation}\label{frob2}
g(s_1,s_2)=s_1s_2-s_1-s_2. 
\end{equation}

However, to find $g(s_1,\dots ,s_n)$, for general $n$, is a difficult problem from the computational point of view;
we refer the reader to \cite{Ram4} for a detailed discussion on the Frobenius number. Klosinski {\it et. al.} used
equation (\ref{frob2}), with particular integers $s_1$ and $s_2$, to show that $(a_1\times a_2)$ can be tiled with $(4\times 6)$ and $(5\times 7)$ if 
$a_1,a_2\ge 2214$.
\vskip .3cm

In this paper, we will use the Frobeniuis number in a more general way to show that a $n$-dimensional rectangle $A$
can be tiled with some set of bricks if the sides of $A$ are larger than a certain function  
(see Theorem \ref{maint}). We use then Theorem \ref{maint} to obtain the following result.

\begin{cor}\label{cor1} Let $a_1,a_2,p,q,r,s\ge 2$ be integers with $r<s$, $g.c.d.(qs,qr,rs,)=1$ and $g.c.d.(p,r)=g.c.d.(p,s)=g.c.d.(r,s)=1$. Then, 
$(a_1\times a_2)$ can be tiled with $(p\times q)$ and $(r\times s)$ if 
$$a_1,a_2\ge \max\{2qrs-(qs+qr+rs)+1,(p-1)(s-1),(r-1)(s-1)\}.$$ 
\end{cor}

In the case when $p=6,q=4,r=5$ and $s=7$, Corollary \ref{cor1} implies that $(a_1\times a_2)$ can be tiled with $(4\times 6)$ and $(5\times 7)$ if 
$a_1,a_2\ge \max\{198,20,24\}=198$, improving the lower bound given in \cite{Klos}. We remark that this lower bound is not optimal. In \cite{Nara}, Narayan and Schwenk 
showed that, in this particular case, it is enough to have $a_1,a_2\ge 33$. However, their tiling constructions allow rotations of both bricks 
(and tilings with more complicated patterns) which is not the case of Corollary \ref{cor1}.

We shall also use Theorem \ref{maint} to prove the following result concerning tilings of squares.
 
\begin{teo}\label{cor2} Let $1<p_1<\cdots <p_{n+1}$ be prime integers. Then, $(\underbrace{a\times \cdots\times a}_n)$ can be tiled with 
$(\underbrace{p_1\times \cdots \times p_1}_n),\dots ,(\underbrace{p_{n+1}\times \cdots\times p_{n+1}}_n)$ if 
$$a>n\prod_{i=1}^{n+1}p_i- \sum\limits_{i=1}^{n+1}{{p_{1}\cdots p_{n+1}}\over {p_i}}\cdot$$
\end{teo}

We finally improve the lower bound given in Theorem \ref{cor2} in some special cases.

\begin{teo}\label{cor3} Let $p>4$ be an odd integer with $3\not | \ p$ and let $a$ be a positive integer. Then, 

$(a\times a)$ can be tiled with $(2\times 2), (3\times 3)$ and $(p\times p)$ if $a\ge 3p+2$. 

Moreover, $(a\times a)$ can be tiled with $(2\times 2), (3\times 3)$ and $(5\times 5)$ if and only if $a\not = 1,7$ and 
with $(2\times 2), (3\times 3)$ and $(7\times 7)$ if and only if $a\not = 1,5,11$.
\end{teo}

A collection of some unpublished work, due to D.A. Klarner, in relation with Theorem \ref{cor3} can be found in \cite{P}. 

\section{Tilings}

We need to introduce some notation and definitions. Let $B=\{b_1,\dots ,b_k\}$ where $b_i$ are positive integers. We
will write $g.c.d.(B)$ instead of $g.c.d.(b_1,\dots ,b_k)$ and $g(B)$ instead of $g(b_1,\dots ,b_k)$.
Let $x_j^i$ be a positive integer for each $j=1,\dots ,n$ and each $i=1,\dots ,n+1$. Let $\{i_1<\dots <i_{k+1}\}\subseteq\{1,\dots ,n+1\}$, $1\le k\le n$.
We define the set

$$X_{k}(i_1,\dots ,i_{k+1})=\left\{{{x_k^{i_1}\cdots x_k^{i_{k+1}}}\over {x_k^{i_{k+1}}}},\dots ,{{x_k^{i_1}\cdots x_k^{i_{k+1}}}\over {x_k^{i_{1}}}}\right\}$$

Note that if $k=n$ then 
$X_n=\{x_n^1\cdots x_n^n,x_n^1\cdots x_n^{n-1}x_n^{n+1},\dots, x_n^2\cdots x_n^{n+1}\}$ is unique.

We denote by $(r;X)$ the rectangle obtained from $X=(x_1\times\cdots\times x_n)$ by sticking together $r$ copies of $X$ along the $n^{th}$-axis, that is,
$(r;X)=(x_1\times\cdots\times x_{n-1}\times rx_n)$. Finally, we denote by $\bar X$ the $(n-1)$-dimensional rectangle obtained from $X=(x_1\times\cdots\times x_n)$
by setting $x_n=0$, that is, ${\bar X}=(x_1\times\cdots\times x_{n-1})$.

\begin{teo}\label{maint} Let $x_j^i\ge 2$ be an integer for each $j=1,\dots ,n$ and each $i=1,\dots ,n+1$, $n\ge 1$. Suppose that 
$g.c.d.(X_k(i_1,\dots ,i_{k+1}))=1$ for any  $\{i_1<\dots <i_{k+1}\}\subseteq\{1,\dots ,n+1\}$, $1\le k\le n$
and let 

$$\hbox{$g_n=\max\limits_{1\le k\le n}\{g(X_k(i_1,\dots ,i_{k+1})) | \{i_1,\dots ,i_{k+1}\}\subseteq\{1,\dots ,n+1\}\}$}$$

Then, $(a_1\times\cdots\times a_n)$ can be tiled with bricks $X^i=(x_1^i\times\cdots\times x_n^i)$, $i=1,\dots ,n+1$ if $a_j> g_n$ for all
$1\le j\le n$.
\end{teo}

{\em Proof.} We shall use induction on $n$. For $n=1$ we have that $g.c.d.(x_1^1,x_1^2)=1$ and thus $g_1=g(x_1^1,x_1^2)=x_1^1x_1^2-x_1^1-x_1^2$. By definition of the Frobenius number, any integer $a_1>g_1$ is of the form $a_1=ux_1^1+vx_1^2$ where $u,v$ are nonnegative integers. Thus, the 1-dimensional rectangle $(a_1)$ (that is, the interval $[0,a_1]$) can be tiled by sticking together $(u;X^1)$ (that is, the interval $[0,ux_1^1]$) and $(v;X^2)$ (that is, the interval $[0,vx_1^2]$). 
\vskip .3cm

We suppose that it is true for $n=m-1\ge 1$ and let $x_j^i$ be a positive 
integer for each $j=1,\dots ,m$ and each $i=1,\dots ,m+1$ with $g.c.d.(X_k(i_1,\dots ,i_{k+1}))=1$ for any 
$\{i_1<\dots <i_{k+1}\}\subseteq\{1,\dots ,m+1\}$, $1\le k\le m$ and let
$X^i=(x_1^i\times\cdots\times x_m^i)$, $i=1,\dots ,m+1$ and $a_j> g_{m+1}$ for all $1\le j\le m$.

By induction $(a_1\times\cdots\times a_{m-1})$ can be tiled with bricks $\bar X^{i_1},\dots ,\bar X^{i_m}$ for any 
$\{i_1<\cdots <i_m\}\subset\{1,\dots ,m+1\}$ since $a_j\ge g_{m-1}$ for all $1\le j\le m-1$. 
\vskip .3cm

We claim that $(a_1\times\cdots\times a_{m-1}\times x_m^{i_1}\cdots x_m^{i_m})$ can be tiled with bricks $X^{i_1},\dots , X^{i_m}$ for any 
$\{i_1<\cdots <i_m\}\subset\{1,\dots ,m+1\}$. 

Indeed, if we consider the rectagle $(a_1\times\cdots\times a_{m-1})$ embedded in $\rere^m$ with $x_m=0$ then
by replacing each brick $\bar X^{i_j}$ used in the tiling of $(a_1\times\cdots\times a_{m-1})$ by
$({{x_m^{i_1}\cdots x_m^{i_m}}\over {x_m^{i_j}}};X^{i_j})$ we obtain a tiling of 
$(a_1\times\cdots\times a_{m-1}\times x_m^{i_1}\cdots x_m^{i_m})$ with bricks $X^{i_1},\dots , X^{i_m}$. 
\vskip .3cm

Now, since $a_m>g_m$ then $a_m=w_{m+1}({x_m^1\cdots x_m^{m+1}\over {x_m^{m+1}}}) +\cdots +w_1({{x_m^1\cdots x_m^{m+1}}\over {x_m^1}})$
where each $w_i$ is a nonnegative integer. By the above claim, $Y^j=(a_1\times\cdots\times a_{m-1}\times {{x_m^{1}\cdots x_m^{m+1}}\over {x_m^j}})$ 
can be tiled with bricks $\{X^{1},\dots ,X^{m+1}\}\setminus X^j$ for each $j=1,\dots ,m+1$. Thus, 
$(a_1\times\cdots\times a_{m-1}\times a_m)$ can be tiled with $X^{1},\dots ,X^{m+1}$ by sticking together bricks $(w_1;Y^1),\dots ,(w_{m+1};Y^{m+1})$ along the
$m^{th}$-axis. \littbox

An old result due to Brauer and Shockley \cite{BS} states that if $s_1,\dots ,s_n$ are positive integers with 
$g.c.d.(s_1,\dots ,s_n)=1$ and $g.c.d.(s_1,\dots ,s_{n-1})=d$ then

\begin{equation}\label{eerr}
g(s_1,\dots ,s_n)=dg\left({{s_1}\over d},\dots ,{{s_{n-1}}\over d},s_n\right)+(d-1)s_n
\end{equation}

\begin{rem}\label{rt} If $s_1,\dots ,s_n$ are positive integers with $g.c.d.(s_1,\dots ,s_n)=1$ and if $s_n$ is representable as a nonnegative 
integer combination of the other $s_i$ then $g(s_1,\dots ,s_n)=g(s_1,\dots ,s_{n-1})$. 
\end{rem}

We shall use equation (\ref{eerr}) and Remark \ref{rt} to prove Corollary \ref{cor1}.

{\em Proof of Corollary \ref{cor1}.} We consider Theorem \ref{maint} with $n=2$, $(x_1^1\times x_2^1)=(p\times q)$, 
$(x_1^2\times x_2^2)=(r\times s)$ and $(x_1^3\times x_2^3)=(s\times r)$. Then,
$(a_1\times a_2)$ can be tiled with $(p\times q)$ and $(r\times s)$ if 

$$a_2,a_2> g_2=\max\{g(p,r),g(p,s),g(r,s), g(qs,qr,sr)\}.$$

By equation (\ref{eerr}), $g(qs,qr,rs)=qg(s,r,rs)+(q-1)rs$ and, by Remark \ref{rt}, 
$g(qs,qr,rs)=qg(s,r)+(q-1)rs=q(sr-s-r)+qrs-rs=2qrs-(qr+qs+rs)$.
The result follows by noticing that $g(p,r)<g(p,s)$ since $r<s$.
\littbox

In order to prove Theorem \ref{maint}, we need the following lemma.

\begin{lem}\label{lls} Let $1<p_1<\cdots <p_{n+1}$ be prime integers, $n\ge 1$. Then,
\begin{enumerate}
\item[(a)] $g\left({{p_{1}\cdots p_{n+1}}\over {p_{n+1}}},\dots, {{p_{1}\cdots p_{n+1}}\over {p_{1}}}\right)=n\prod\limits_{i=1}^{n+1}p_i- \sum\limits_{i=1}^{n+1}{{p_{1}\cdots p_{n+1}}\over {p_i}}\cdot$

\item[(b)] For each $k=1,\dots ,n$,

$$g\left({{p_{n-k+1}\cdots p_{n+1}}\over {p_{n+1}}},\dots, {{p_{n-k+1}\cdots p_{n+1}}\over {p_{n-k+1}}}\right)\ge g\left({{p_{i_1}\cdots p_{i_{k+1}}}\over {p_{i_{k+1}}}},\dots, 
{{p_{i_1}\cdots p_{i_{k+1}}}\over {p_{i_{1}}}}\right)$$

for all $\{i_1<\dots <i_{k+1}\}\subseteq \{1,\dots ,n+1\}$ with equality if and only if $\{i_1,\dots ,i_{n+1}\}= \{1,\dots ,n+1\}$.

\end{enumerate}
\end{lem}

\proof (a) By induction on $n$. For $n=1$ we have that $g(p_1,p_2)=p_1p_2-p_1-p_2$ finding equation (\ref{frob2}). We suppose that it is true for $n=m$. 
We assume that $n=m+1$, by equation (\ref{eerr}), we have

$$g\left({{p_{1}\cdots p_{m+1}}\over {p_{m+1}}},\dots, {{p_{1}\cdots p_{m+1}}\over {p_{1}}}\right)=p_1g\left({{p_{2}\cdots p_{m+1}}\over {p_{m+1}}},\dots ,{{p_{2}\cdots p_{m+1}}\over {p_{2}}}, {{p_{1}\cdots p_{m+1}}\over {p_{1}}}\right)+(p_1-1){{p_{1}\cdots p_{m+1}}\over {p_{1}}}$$

Since ${{p_{1}\cdots p_{m+1}}\over {p_{1}}}=p_2{{p_{2}\cdots p_{m+1}}\over {p_{2}}}$ then, by Remark \ref{rt},

$$g\left({{p_{2}\cdots p_{m+1}}\over {p_{m+1}}},\dots ,{{p_{2}\cdots p_{m+1}}\over {p_{2}}}, {{p_{1}\cdots p_{m+1}}\over {p_{1}}}\right)=g\left({{p_{2}\cdots p_{m+1}}\over {p_{m+1}}},\dots, {{p_{2}\cdots p_{m+1}}\over {p_{1}}}\right)$$
 and thus

$$\begin{array}{ll}
g\left({{p_{1}\cdots p_{m+1}}\over {p_{m+1}}},\dots, {{p_{1}\cdots p_{m+1}}\over {p_{1}}}\right) & = p_1g\left({{p_{2}\cdots p_{m+1}}\over {p_{m+1}}},\dots, {{p_{2}\cdots p_{m+1}}\over {p_{1}}}\right)+(p_1\cdots p_{m+1})-{{p_{1}\cdots p_{m+1}}\over {p_{1}}}\\
& = (m-1)\prod\limits_{i=1}^{m+1}p_i- \sum\limits_{i=2}^{m+1}{{p_{1}\cdots p_{m+1}}\over {p_i}}+\prod\limits_{i=1}^{m+1}p_i-{{p_{1}\cdots p_{m+1}}\over {p_{1}}}\\
& = p_1\left((m-1)\prod\limits_{i=2}^{m+1}p_i- \sum\limits_{i=2}^{m+1}{{p_{2}\cdots p_{m+1}}\over {p_i}}\right)+\prod\limits_{i=1}^{m+1}p_i-{{p_{1}\cdots p_{m+1}}\over {p_{1}}}\\
& = m\prod\limits_{i=1}^{m+1}p_i- \sum\limits_{i=1}^{m+1}{{p_{1}\cdots p_{m+1}}\over {p_i}}\\
\end{array}$$

(b) By induction on $k$. For $k=1$ we have $g(p_{i_1},p_{i_2})=(p_{i_1}-1)(p_{i_2}-1)\le (p_{n}-1)(p_{n+1}-1)=g(p_n,p_{n+1})$ since $p_{i_1}\le p_{n}$  and 
$p_{i_2}\le p_{n+1}$ for all $\{i_1<i_2\}\subset\{1,\dots ,n+1\}$. We suppose that it is true for $k=m$ and assume that $k=m+1\le n$. By equation (\ref{eerr}) we have

$$\begin{array}{ll}
g\left({{p_{n-m}\cdots p_{n+1}}\over {p_{n+1}}},\dots, {{p_{n-m}\cdots p_{n+1}}\over {p_{n-m}}}\right)& = p_{n-m}g\left({{p_{n-m+1}\cdots p_{n+1}}\over {p_{n+1}}},\dots, {{p_{n-m+1}\cdots p_{n+1}}\over {p_{n-m-1}}}\right)\\
&+(p_{n-m}-1){{p_{n-m}\cdots p_{n+1}}\over {p_{n-m}}}\\
\end{array}$$

and

$$\begin{array}{ll}
g\left({{p_{i_{1}}\cdots p_{i_{m+2}}}\over {p_{i_{m+2}}}},\dots, {{p_{i_{1}}\cdots p_{i_{m+2}}}\over {p_{i_{1}}}}\right)& =p_{i_{1}}g\left({{p_{i_{2}}\cdots p_{i_{m+2}}}\over {p_{i_{m+2}}}},\dots, {{p_{i_{2}}\cdots p_{i_{m+2}}}\over {p_{i_{2}}}}\right)\\
&+(p_{i_{1}}-1){{p_{i_{1}}\cdots p_{i_{m+2}}}\over {p_{i_{1}}}}\\
\end{array}$$

By induction, we obtain
  
$$g\left({{p_{n-m+1}\cdots p_{n+1}}\over {p_{n+1}}},\dots, {{p_{n-m+1}\cdots p_{n+1}}\over {p_{n-m-1}}}\right)\ge g\left({{p_{i_{2}}\cdots p_{i_{m+2}}}\over {p_{i_{m+2}}}},\dots, {{p_{i_{2}}\cdots p_{i_{m+2}}}\over {p_{i_{2}}}}\right)$$

Moreover, since $i_1<\cdots <i_{m+2}\le n+1$ then $i_{m+2-j}\le n+1-j$ for each $j=0,\dots ,m+1$ and thus $p_{i_{m+1-j}}\le p_{n-j}$. Therefore,
$p_{i_1}=p_{i_{m+2-(m+1)}}\le p_{n+1-(m+1)}=p_{n-m}$ and

$$\begin{array}{ll}
{{p_{n-m}\cdots p_{n+1}}\over {p_{n-m}}} & =p_{n-m+1}p_{n-m+2}\cdots p_{n+1}\\
& = p_{n+1-m}p_{n+1-(m-1)}\cdots p_{n+1}\\
&\ge p_{i_{m+2-m}}p_{i_{m+2-(m-1))}}\cdots p_{i_{m+2}}\\
& =p_{i_{2}}p_{i_{3}}\cdots p_{i_{m+2}}\\
&= {{p_{i_1}\cdots p_{i_{m+2}}}\over {p_{i_{1}}}}
\end{array}$$

Obtaining that 
$$g\left({{p_{n-m}\cdots p_{n+1}}\over {p_{n+1}}},\dots, {{p_{n-m}\cdots p_{n+1}}\over {p_{n-m}}}\right)\ge g\left({{p_{i_{1}}\cdots p_{i_{m+2}}}\over {p_{i_{m+2}}}},\dots, {{p_{i_{1}}\cdots p_{i_{m+2}}}\over {p_{i_{1}}}}\right).$$
\littbox

We may now prove Theorem \ref{cor2}.

{\em Proof of Theorem \ref{cor2}.}  Let $1<p_1<\cdots <p_{n+1}$ be prime integers, $n\ge 1$. We consider Theorem \ref{maint}, where, for each $i=1,\dots ,n+1$, 
we let $x_i^j=p_i$ for all $1\le j\le n$. Then, $(\underbrace{a\times \cdots\times a}_n)$ can be tiled with 
$(\underbrace{p_1\times \cdots \times p_1}_n),\dots ,(\underbrace{p_{n+1}\times \cdots\times p_{n+1}}_n)$ if 
$$\hbox{$a>g_n=\max\limits_{1\le k\le n}\{g(X_k(i_1,\dots ,i_{k+1})) | \{i_1,\dots ,i_{k+1}\}\subseteq\{1,\dots ,n+1\}$.}$$

By Lemma \ref{lls} (b), for each fixed $k=1,\dots ,n$

$$\begin{array}{ll}
\bar g_k&=\max\{g(X_k(i_1,\dots ,i_{k+1})) | \{i_1,\dots ,i_{k+1}\}\subseteq\{1,\dots ,n+1\}\}\\
&=g\left({{p_{n-k+1}\cdots p_{n+1}}\over {p_{n+1}}},\dots, {{p_{n-k+1}\cdots p_{n+1}}\over {p_{n-k+1}}}\right)
\end{array}$$

Moreover, for each $k=1,\dots ,n$, $\bar g_{k+1}>\bar g_k$ since

$$\begin{array}{ll}
\bar g_{k+1}=g\left({{p_{n-k}\cdots p_{n+1}}\over {p_{n+1}}},\dots, {{p_{n-k}\cdots p_{n+1}}\over {p_{n-k}}}\right)& = 
p_{n-k}g\left({{p_{n-k+1}\cdots p_{n+1}}\over {p_{n+1}}},\dots, {{p_{n-k+1}\cdots p_{n+1}}\over {p_{n-k-1}}}\right)\\
&+(p_{n-k}-1){{p_{n-k}\cdots p_{n+1}}\over {p_{n-k}}}\\
&> g\left({{p_{n-k+1}\cdots p_{n+1}}\over {p_{n+1}}},\dots, {{p_{n-k+1}\cdots p_{n+1}}\over {p_{n-k-1}}}\right)=\bar g_k\\
\end{array}$$

So, $g_{n}=\max\limits_{k=1,\dots ,n} \{\bar g_k\}$ and, by Lemma \ref{lls} (a), 

$$g_{n}=g\left({{p_{1}\cdots p_{n+1}}\over {p_{n+1}}},\dots, {{p_{1}\cdots p_{n+1}}\over {p_{1}}}\right)=n\prod\limits_{i=1}^{n+1}p_i- \sum\limits_{i=1}^{n+1}{{p_{1}\cdots p_{n+1}}\over {p_i}}\cdot$$
\littbox

\begin{pro}\label{llls} Let $L,a,b,c$ and $r$ be positive integers with $b|r$ and such that $r=x_1a+x_2c$ for some integers $x_1,x_2\ge 0$ 
and $Lc=y_1a+y_2b$ for some integers $y_1,y_2\ge 0$. Then, 
$(r+ac\times r+ac)$ and $(Lc+kab\times Lc+kab)$ can be tiled with $(a\times a),(b\times b)$ and $(c\times c)$ for any integer $k\ge 1$.
\end{pro}

\proof Suppose that $b|r$. By Theorem \ref{kler}, we have that 

$\bullet$ $(r\times r)$ can be tiled with $(b\times b)$

$\bullet$ $(ac\times ac)$ can be tiled with $(a\times a)$,

$\bullet$ $(ac\times r)$ can be tiled with $(a\times a)$ and $(c\times c)$, 

$\bullet$ $(Lc\times Lc)$ can be tiled with $(c\times c)$,

$\bullet$ $(Lc\times kab)$ can be tiled with $(a\times a)$ and $(b\times b)$ and

$\bullet$ $(kab\times kab)$ can be tiled with $(a\times a)$ (or $(b\times b)$).

\begin{figure}[h]
\begin{center}
\leavevmode
\epsfxsize=270pt
{\epsffile{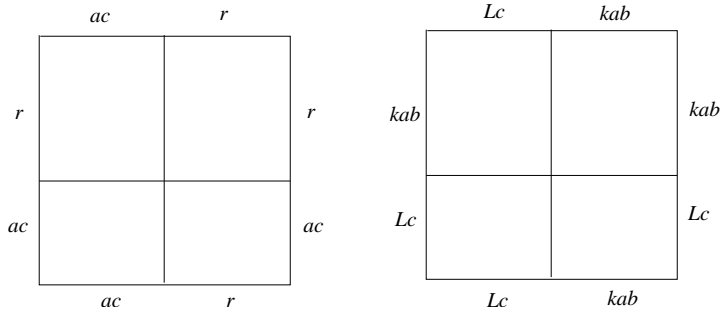}}
\caption{Compositions of tilings}\label{fig1}
\end{center}
\end{figure}

The results folllow by sticking together copies of the above rectangles as shown in Figure \ref{fig1}.\littbox

{\em Proof of Theorem \ref{cor3}.} By Theorem \ref{kler}, we have that $(a\times a)$ can be tiled with $(2\times 2)$ and $(3\times 3)$ if
$a\equiv 0 \bmod 2$ or $a\equiv 0 \bmod 3$. So, we only need to show that 
$(a\times a)$ can be tiled with $(2\times 2),(3\times 3)$ and $(p\times p)$ for any odd integer $a\ge 3p+2$ with $a\equiv 1$ or $2 \bmod 3$.
We have two cases.
\vskip .3cm

Case A) $p\equiv 1\bmod 3$. Let $s=p-1+k3\ge p+2$ for any integer $k\ge 1$. Since $s> g(2,p)=p-1$ then there exist nonnegative integers $u$ and $v$ 
such that $s=2u+pv$. Since $3|s$ then, by Proposition \ref{llls}, $(a\times a)=(s+2p\times s+2p)$ can be tiled with  $(2\times 2),(3\times 3)$ and $(p\times p)$
for any $a=s+2p\ge p+2+2p=3p+2$ and $a\equiv 2 \bmod 3$.

Now, since $p=3k+1$ for some integer $k\ge 1$ then for $p>3$ we have that $p= (k-1)3+2(2)$ and 
by Proposition \ref{llls} (with $L=1$), we have that $(a\times a)=(p+k6\times p+k6)$
can be tiled with $(2\times 2),(3\times 3)$ and $(p\times p)$ for any odd integer $a=p+k6\ge p+6$ with $a\equiv 1\bmod 3$.
\vskip .3cm

Case B) $p\equiv 2\bmod 3$. Let $s=p-2+k3\ge p+1$ for any integer $k\ge 1$. Since $s> g(2,p)=p-1$ then there exist nonnegative integers $u$ and $v$ 
such that $s=2u+pv$. Since $3|s$ then, by Proposition \ref{llls}, $(a\times a)=(s+2p\times s+2p)$ can be tiled with  $(2\times 2),(3\times 3)$ and $(p\times p)$
for any $a=s+2p\ge 2p+1+p=3p+1$ and $a\equiv 1 \bmod 3$.

Now, since $p=3k+2$ for some integer $k\ge 1$ then for $p>3$ we have that $p= (k-1)3+2(2)$ and 
by Proposition \ref{llls} (with $L=1$), we have that $(a\times a)=(p+k6\times p+k6)$
can be tiled with $(2\times 2),(3\times 3)$ and $(p\times p)$ for any odd integer $a=p+k6\ge p+6$ with $a\equiv 2\bmod 3$.
\vskip .3cm

\begin{figure}[h]
\begin{center}
\leavevmode
\epsfxsize=160pt
{\epsffile{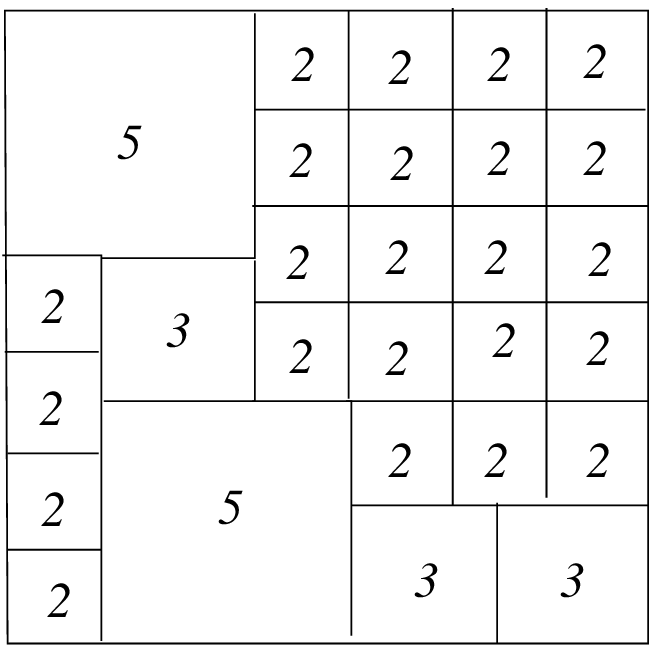}}
\caption{Tiling of $(13\times 13)$ with bricks $(2\times 2),(3\times 3)$ and $(5\times 5)$}\label{fig2}
\end{center}
\end{figure}

Let us set $p=5$. It is clear that $(1\times 1)$ and $(7\times 7)$ cannot be tiled with $(2\times 2),(3\times 3)$ and $(5\times 5)$. 
By the above cases, we have that $(a\times a)$ can be tiled with $(2\times 2),(3\times 3)$ and $(5\times 5)$ if $a\ge 3p+2=17$
and, by Theorem \ref{kler}, $(a\times a)$ can be tiled with $(2\times 2)$ and $(3\times 3)$ if $a\equiv 0 \bmod 2$ or $a\equiv 0 \bmod 3$. 
This leave us the cases when $a=5,11$ and $13$.  If $a=5$ is trivial.
$(11\times 11)$ can be tiled with $(2\times 2),(3\times 3)$ and $(5\times 5)$ since, by tha above case B, the result is true for 
any odd integer $a\ge p+6=11$ and $a\equiv 2 \bmod 3$. Finally,
$(13\times 13)$ can be tiled as it is illustrated in Figure \ref{fig2}.

\begin{figure}[h]
\begin{center}
\leavevmode
\epsfxsize=170pt
{\epsffile{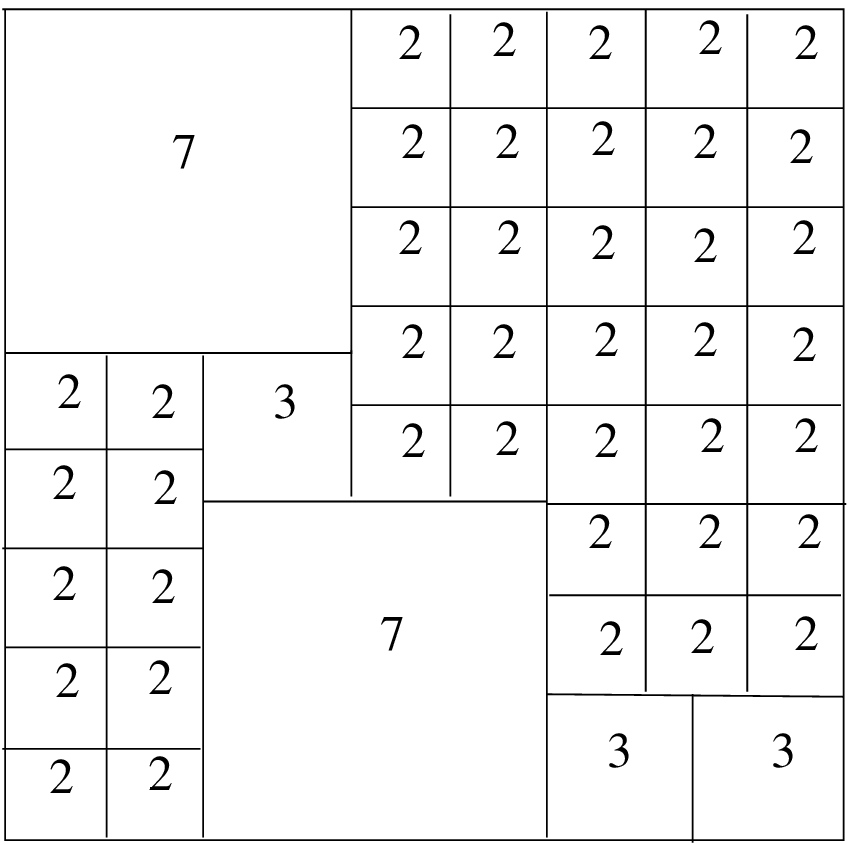}}
\caption{Tiling of $(17\times 17)$ with bricks $(2\times 2),(3\times 3)$ and $(7\times 7)$ }\label{fig3}
\end{center}
\end{figure}

Let us set $p=7$. It is clear that $(1\times 1)$, $(5\times 5)$ and $(11\times 11)$ cannot be tiled with $(2\times 2),(3\times 3)$ and $(7\times 7)$. 
By the above cases, we have that $(a\times a)$ can be tiled with $(2\times 2),(3\times 3)$ and $(7\times 7)$ if $a\ge 3p+2= 23$
and, by Theorem \ref{kler}, $(a\times a)$ can be tiled with $(2\times 2)$ and $(3\times 3)$ if $a\equiv 0 \bmod 2$ or $a\equiv 0 \bmod 3$. 
This leave us the cases when $a=7,13,17$ and $19$. If $a=7$ is trivial. 
$(13\times 13)$ and $(19\times 19)$ both can be tiled since, by the above case A, the result is true for any odd integer $a\ge p+6=13$ with $a\equiv 1\bmod 3$. 
Finally, $(17\times 17)$ can be tiled as it is illustrated in Figure \ref{fig3}.

\littbox

{\bf Acknowledgement}

I would like to thank Y. O. Hamidoune for several helpful conversations.


\begin{thebibliography}{99}

\bibitem{Bar} F.W. Barnes, Algebraic theory of bricks packing II, {\em Discrete Mathematics} {\bf 42} (1982), 129-144.

\bibitem{BS} A. Brauer and J.E. Shockley, On a problem of Frobenius, {\em J. Reine Angewandte Math.} {\bf 211} (1962), 215-220.

\bibitem{BW} R.J. Brower and T.S. Michael, When can you tile a box with translates of two given rectangular bricks?, 
{\it Elec. J. of Comb.} {\bf 11} (2004), \# N7.

\bibitem{deB} N.G. de Bruijn, Filling boxes with bricks, {\it Amer. Math. Monthly} {\bf 76} (1969), 37-40.

\bibitem{Fri} J. Fricke, Quadratzerlegung eines Rechtecks, {\it Math. Semesterber.} {\bf 42} (1995), 53-62.

\bibitem{J} S.M. Johnson, A linear diophantine problem, {\it Can. J. Math.}
{\bf 12} (1960), 390-398.

\bibitem{Kla} D.A. Klarner, Packing a rectangle with congruent $n$-ominoes, {\it J. Comb. Theory} {\bf 7} (1969), 107-115.

\bibitem{Klos} L.F. Klosinski, G.L. Alexanderson and L.C. Larson, The Fifty-Second William Lowell Putman Mathematical Competition,
{\it Amer. Math. Month.} {\bf 9} (1992), 715-724.

\bibitem{Nara} D. A. Narayan and A.J. Schwenk, Tiling large rectangles
{\it Mathematics Magazine} {\bf 75}(5) (2002), 372-380.

\bibitem{P} T. Plambeck, Web site http://www.plambeck.org/oldhtml/mathematics/klarner/

\bibitem{Ram4} J.L. Ram\'{\i}rez Alfons\'{\i}n, The Diophantine Frobenius Problem, 
{\it Oxford Lectures Series in Mathematics and its Applications} {\bf 30}, Oxford University Press, (2005). 

\end{thebibliography}
\end{document}